December 07, 2016

# A note on the higher derivatives of the function 1/(exp(x) - 1)

Khristo N. Boyadzhiev

*Department of Mathematics and Statistics, Ohio Northern University, Ada, OH 45810, USA*
k-boyadzhiev@onu.edu

**Abstract**. In this brief note we compare two formulas for the higher order derivatives of the function 1/(exp(x) -1). We also provide an integral representation for these derivatives and obtain a classical formula relating zeta values and Bernoulli numbers.

## 1. Two formulas

Let $S(n,k)$ be the Stirling numbers of the second kind (see [9]). The following formula has been proved by induction in several publications (for instance, [1], [10], [11], [12], [13]). For $n \geq 1$:

$$(1) \qquad \left(\frac{d}{dx}\right)^n \frac{1}{e^x - 1} = (-1)^n \sum_{k=1}^{n+1} S(n+1,k)(k-1)! \left(\frac{1}{e^x - 1}\right)^k .$$

Another formula for the same derivative is this

$$(2) \qquad \left(\frac{d}{dx}\right)^n \frac{1}{e^x - 1} = (-1)^n \frac{e^x}{e^x - 1} \sum_{k=1}^{n} S(n,k)\, k! \left(\frac{1}{e^x - 1}\right)^k .$$

Formula (2) turns into (1) immediately when we write

$$\frac{e^x}{e^x - 1} = \frac{e^x - 1 + 1}{e^x - 1} = 1 + \frac{1}{e^x - 1} ,$$

then split the sum in two sums, and use the well-known property of the Stirling numbers

$$S(n+1,k) = k\, S(n,k) + S(n,k-1) .$$

## 2. Proof of the second formula and discussion

It was shown in [3] and [8] that for $n \geq 1$, $x \neq 1$,

$$(3) \qquad \left(\frac{d}{dx}\right)^n \frac{1}{1-x} = \frac{1}{1-x} \sum_{k=1}^{n} S(n,k)\, k! \left(\frac{x}{x-1}\right)^k = \frac{1}{1-x} \omega_n \left(\frac{x}{x-1}\right)$$

where the polynomials $\omega_n(x) = \sum_{k=0}^{n} S(n,k)\, k!\, x^k$ were called in [3] geometric polynomials. The geometric polynomials and the generalized geometric polynomials $\omega_{n,p}$ introduced in [3] appear in many applications (see, for example [5], [6] and the references in [6]).

Note that equation (3) for $|x| < 1$ can be written also in the form (cf. (7.46) in [9])

$$\left(\frac{d}{dx}\right)^n \frac{1}{1-x} = \sum_{k=0}^{\infty} m^n x^m = \frac{1}{1-x} \omega_n \left(\frac{x}{x-1}\right)$$

Formula (2) follows immediately from (3) when we replace $x$ by $e^{-x}$

$$\left(\frac{d}{dx}\right)^n \frac{1}{e^x - 1} = \left(\frac{d}{dx}\right)^n \frac{e^{-x}}{1 - e^{-x}} = \left(\frac{d}{dx}\right)^n \left(-1 + \frac{1}{1 - e^{-x}}\right)$$

$$= (-1)^n \frac{1}{1 - e^{-x}} \sum_{k=1}^{n} S(n,k)\, k! \left(\frac{e^{-x}}{1 - e^{-x}}\right)^k .$$

This formula was used in [4] for computing the derivative polynomials for $\tan x$, $\cot x$, $\sec x$, $\tanh x$, $\operatorname{sech} x$, and $\coth x$.

The geometric polynomials $\omega_n$ were used in [7] to obtain the Maclaurin series in $t$

$$\frac{1}{\mu e^{\lambda t} + 1} = \frac{1}{\mu + 1} \sum_{n=0}^{\infty} \lambda^n \omega_n \left(\frac{-\mu}{\mu + 1}\right) \frac{t^n}{n!}$$

with parameters $\lambda$, $\mu$ .

## 3. Integral representation and Bernoulli numbers

**Proposition**. *For every nonnegative integer $n$ and every $x > 0$ we have*

(4)
$$\sum_{k=1}^{\infty} k^n e^{-kx} = \frac{n!}{x^{n+1}} + 2(-1)^n \int_0^{\infty} \frac{t^n}{e^{2\pi t} - 1} \sin\left(xt + \frac{n\pi}{2}\right) dt \, .$$

Equation (4) is a known result that can be found, for example, in the works of Ramanujan (see Entry 2 on p. 411 in [2]). A simple proof of (4) can be obtained this way:

Consider the equation (Fourier sine transform formula)

$$\int_0^{\infty} \frac{\sin(xt)}{e^{2\pi t} - 1} dt = \frac{1}{4}\coth\frac{x}{2} - \frac{1}{2x} = \frac{1}{4}\left(1 + 2\frac{e^{-x}}{1 - e^{-x}}\right) - \frac{1}{2x}$$

which, by using geometric series, can be written in the form

$$\sum_{k=1}^{\infty} e^{-kx} - \frac{1}{x} = 2\int_0^{\infty} \frac{\sin(xt)}{e^{2\pi t} - 1} dt - \frac{1}{2} \, .$$

Differentiating this equation $n$ times for $x$ yields (4).

Now, from equation (4) and the equation

$$\left(\frac{d}{dx}\right)^n \frac{1}{e^x - 1} = \left(\frac{d}{dx}\right)^n \frac{e^{-x}}{1 - e^{-x}} = \left(\frac{d}{dx}\right)^n \sum_{k=1}^{\infty} e^{-kx} = (-1)^n \sum_{k=1}^{\infty} k^n e^{-kx}$$

we derive the integral representation

(5)
$$\left(\frac{d}{dx}\right)^n \frac{1}{e^x - 1} = \frac{(-1)^n n!}{x^{n+1}} + 2\int_0^{\infty} \frac{t^n}{e^{2\pi t} - 1} \sin\left(xt + \frac{n\pi}{2}\right) dt \, ,$$

or

(6)
$$\left(\frac{d}{dx}\right)^n \left(\frac{1}{e^x - 1} - \frac{1}{x}\right) = 2\int_0^{\infty} \frac{t^n}{e^{2\pi t} - 1} \sin\left(xt + \frac{n\pi}{2}\right) dt$$

From (5) we compute the interesting limit ($n \geq 1$)

(7)
$$\lim_{x \to 0} \left(\frac{d}{dx}\right)^n \left(\frac{1}{e^x - 1} - \frac{1}{x}\right) = \frac{2n!}{(2\pi)^{n+1}} \sin\left(\frac{n\pi}{2}\right) \zeta(n+1)$$

where

$$\zeta(s) = \frac{1}{\Gamma(s)} \int_0^{\infty} \frac{t^{s-1}}{e^t - 1} dt$$

is the Riemann zeta function.

It is natural to connect this result to the Bernoulli numbers defined by the generating function

$$\frac{x}{e^x - 1} = \sum_{n=0}^{\infty} \frac{B_n}{n!} x^n .$$

Since $B_0 = 1$, we can write this expansion in the form

$$\frac{1}{e^x - 1} - \frac{1}{x} = \sum_{n=0}^{\infty} \frac{B_{n+1}}{(n+1)!} x^n$$

and (7) implies

$$\frac{B_{n+1}}{(n+1)!} = \frac{2}{(2\pi)^{n+1}} \sin\left(\frac{n\pi}{2}\right) \zeta(n+1) .$$

When $n+1$ is odd both sides are zeros. When $n+1 = 2m$ this turns into a classical result of Euler

$$B_{2m} = \frac{(-1)^{m+1} 2(2m)!}{(2\pi)^{2m}} \zeta(2m) .$$